\documentclass{amsart}
\usepackage{amssymb}
\usepackage{graphics}

\theoremstyle{plain}
\newtheorem{lemma}{Lemma}

\newtheorem{theorem}{Theorem}
\theoremstyle{definition}
\newtheorem{definition}{Definition}
\newtheorem{corollary}{Corollary}
\theoremstyle{definition}
\newtheorem{remark}{Remark}

\begin{document}
\title[Circle actions on symplectic manifolds]{On the construction of certain 6-dimensional symplectic manifolds with  Hamiltonian circle actions}
\author{Hui Li}
\address{Mathematics Department \\
         University of Illinois \\
         Urbana-Champaign, IL 61801 }
\email{hli@math.uiuc.edu}
%\thanks{The author would like to thank...}
\begin{abstract}
   Let $(M, \omega)$ be a connected, compact 6-dimensional symplectic manifold equipped with a semi-free Hamiltonian $S^1$ action such that
   the fixed point set consists of isolated points or surfaces. Assume dim $H^2(M)<3$, in \cite{L}, we defined a certain invariant of such
  spaces which consists of fixed point data and twist type, and we divided the possible values of these invariants into six ``types''.
   In this paper, we construct such manifolds with these ``types''. As a consequence, we have a precise list of the values of these invariants.  
\end{abstract}
 \maketitle
 \section{introduction}
 Assume $(M, \omega)$ is a connected, compact 6-dimensional symplectic manifold equipped with a semi-free (free outside fixed point set)
 Hamiltonian $S^1$ action such that
   the fixed point set consists of isolated points or surfaces. For the case dim $H^2(M)<3$, in \cite{L}, we defined a certain invariant which
  consists of the ``fixed point data'' and ``twist type'', and we gave  a complete list
   of the possible values of these invariants. We have seen that only six ``types'' can occur (see Theorem~\ref{thm1} below). We showed the existences for a few cases
   by giving  explicit examples, such as coadjoint orbits or toric varieties. In this paper, we illustrate a method for constructing
  such manifolds. As a consequence, this proves the existences of the manifolds with those prescribed types as in \cite{L}.
  The techniques we use in this paper are quite different from those we used in \cite{L}. We mainly 
   use  symplectic topological method and related results on symplectic 4-manifolds. 
   
     Recall that the moment map $\phi$ is a perfect Morse-Bott function. The fixed point sets of the circle action are exactly the critical
   sets of $\phi$. Each  critical set has even index. $\phi$ has a unique local minimum and a unique local maximum (\cite{A}).
   Recall the following definitions from \cite{L}:
  \begin{definition}
  Let $(M, \omega)$ be a connected, compact 6-dimensional symplectic manifold equipped with a semi-free Hamiltonian $S^1$ action such that
   the fixed point set consists of isolated points or surfaces. Assume dim $H^2(M)<3$.
  For each such manifold, we  define the $\bold{fixed\, point\, data}$ to be the diffeomorphism type of the fixed point sets, their indices,
 and $b_{min}$, where $b_{min}$ denotes the 1st Chern number of the normal bundle of the minimum when the minimum is a surface.
\end{definition}
         
    If the minimum (or the maximum) is a surface, then for a regular value right above the minimum (or right below the maximum), the
   reduced space is diffeomorphic to the total space of an $S^2$ fibration over the minimum (or over the maximum). 
 When the fixed point set consists only of surfaces, 
   the reduced spaces at the non-extremal values are all diffeomorphic, and 
  there is a canonical isomorphism between their (co)homology groups (see \cite{GS}).
    \begin{definition}
   Let $(M,\omega)$ be a compact connected 6-dimensional symplectic manifold equipped with a semi-free Hamiltonian circle action.
  Assume the fixed point set consists only of surfaces. If the fibers of the reduced spaces at a regular level right
 above the minimum and right below the maximum can be represented by the same cohomology class of the reduced space,
 we say there is $\bold{no\, twist}$; otherwise, we say there is a $\bold{twist}$.
 \end{definition}

 In \cite{L}, we have shown the following 
 
    \begin{theorem}\label{thm1}
     Let $(M,\omega)$ be a  connected compact 6 dimensional symplectic manifold equipped with a semi-free Hamiltonian circle action 
     such that the fixed point set consists of isolated  points or  surfaces.  Assume $\dim$ $H^2(M)<3$. 
     Then the possible fixed point data and the twist type are:

    \begin{figure}[h!]
    \scalebox{.60}{\includegraphics{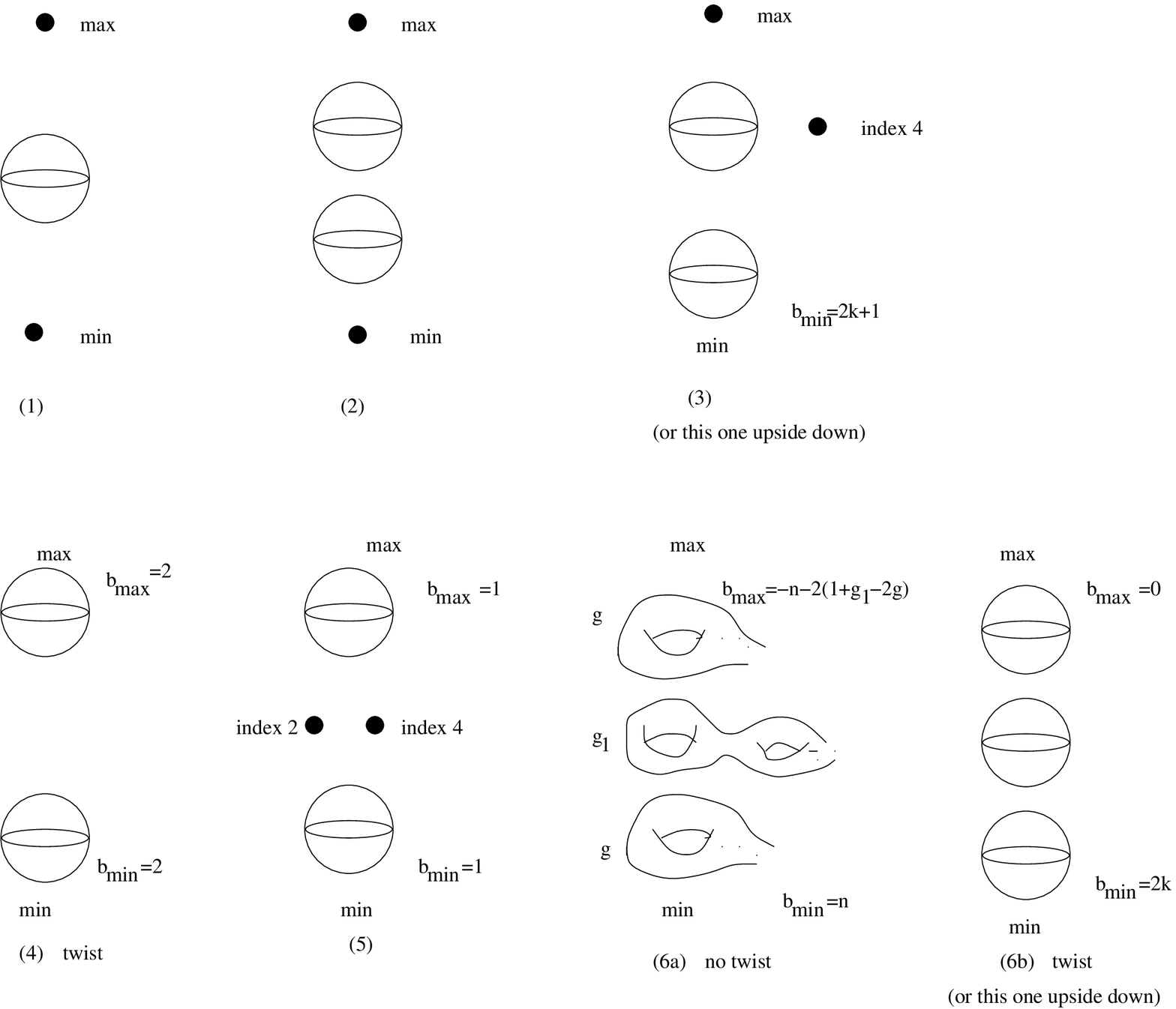}}
    \end{figure}

 In the picture, $g$, $g_1$ are the genuses of the surfaces, and $b_{max}$ is the 1st Chern number of the normal bundle of the  maximum.
    \end{theorem}

\begin{remark}\label{rem1}
 The picture gives  us an order of how the moment map crosses the fixed point sets. For instance, manifold of type (1) means that the manifold has an isolated minimum,
an index 2 sphere, and an isolated maximum.

    For manifold of type (2), the two index 2 spheres can not be on the same level of the moment map.

    For manifolds of type (3) and (5), the moment map can cross the two non-extremal fixed point sets  in any order or the 2 fixed point sets are on the same level of the moment map.
\end{remark}

 \begin{remark}\label{rem1'}
     In \cite{L}, we have seen that when there is an index 2 surface or the maximum is a surface in a certain ``type'', 
  the 1st Chern numbers of the normal bundles of these surfaces are uniquely determined by the fixed point data and the twist type. (This is the reason that
  $b_{max}$ is not included in the definition of the $\bold{fixed\, point\, data}$.)  For instance, for manifolds
 of type (1), the 1st Chern numbers of the positive normal bundle and the negative normal bundle of the index 2 sphere are both $2$; for manifolds of
  type (2), the 1st Chern numbers of the positive normal bundle and the negative normal bundle of the 1st (counting from bottom) index 2 sphere are $0$ and
 $1$ respectively, and those of the second index 2 sphere are respectively $1$ and $0$.    
  Notice that in type (4) and type (5), $b_{min}$ and $b_{max}$ have fixed values, in types (3), (6a) and (6b), on the other hand,  $b_{min}$ can vary, and
 in particular, for type (6a), the genuses $g$ and $g_1$ can vary. Different values of these ``variables'' of course give different manifolds.   
\end{remark}

  \begin{remark}\label{rem2}
      We have shown  in \cite{L} that there exists a manifold of type (1) which is a coadjoint orbit of $SO(5)$.
  We constructed  toric varieties for  type (3) when $b_{min}=1$ and $b_{min}=3$, for type (4) and for  type (6b) with $b_{min}=2$.
  Moreover, we have shown that all manifolds of type (4) are diffeomorphic to $\mathbb{C}P^3$, and all manifolds of type (6b)  with $b_{min}=2$
  are diffeomorphic to the corresponding toric variety which is  $\mathbb{C}P^3$ blow up at a point.
     \end{remark}

     In  \cite{L} and in this paper, the $\bold{fixed\, point\, data}$
   has a certain meaning, for the sake of convenience. More generally, similar data (looser or stronger) is an important tool for the studying
  of symplectic Lie group actions, especially for the studying of Hamiltonian Lie group actions. It captures important geometrical and topological
  information about the manifolds. 
    As we have seen in  \cite{L}, for a fixed ``type'' in Theorem~\ref{thm1}, the $\bold{fixed\, point\, data}$ uniquely determines the equivariant (ordinary) cohomology ring 
  and equivariant (ordinary) Chern classes of the manifold. More significantly, for some cases, this data classifies
   the manifolds up to diffeomorphism (see  \cite{L}). Recall that for all the types of manifolds in Theorem~\ref{thm1},
 the reduced spaces (except at the extremal values) are $\mathbb{C}P^2$ or $\bold{ruled\, manifolds}$. We understand well the symplectic
 structures of these 4-manifolds and their symplectic submanifolds.  In this paper, we use this special feature to construct
 manifolds with the types prescribed in Theorem~\ref{thm1}. We will show that for each ``type'' in Theorem~\ref{thm1}, except for type
 (6a) with $1+g_1-2g\leq 0$,
  this ``local'' data (related in a global way, however) allows us to construct 
   a ``global'' manifold  with the prescribed data. In fact, except for type
 (6a) with $1+g_1-2g\leq 0$, this data uniquely determines the inverse image (under the moment map) of a small neighborhood
   of each critical or regular value of the moment map up to equivariant symplectomorphism (see Section 2). Let us call these
 ``local pieces''. A ``local piece'' for the minimum or the maximum can always be constructed. The main problem is to construct
 a ``local piece'' for an index 2 surface. In this case, we need 4 pieces of information $(M_0, \omega_0, P, X)$ (see Section 2 for its meaning),
  where $M_0$ is the reduced space at this critical level, and $X$ has to be a symplectic submanifold of $M_0$. More precisely, $X$ is
 the symplectically embedded image of the index 2 surface in  $M_0$ (see Section 2). The prescribed type must make this to be true. Of course, 
 the prescribed type must guarantee that each reduced space is symplectic, i.e., the cohomology classes of the reduced symplectic
 forms (we can only know the cohomology classes of the forms) do contain symplectic forms (see Section 3).  For (6a), if
 $1+g_1-2g\leq 0$, the dual class (which can be computed from the ``type'') of the surface in the reduced space can not represent the symplectically
 embedded image of the index 2 surface. 
 That is the reason we failed to construct such  manifolds for type (6a) when  $1+g_1-2g\leq 0$. 
  In order to construct the manifolds globally,
  we need to perform some gluing along some regular levels. The main reason we can perform the gluing is that for a regular level of 
  the moment map,
  two different symplectic forms in the same cohomology class of the reduced space (which is $\mathbb{C}P^2$ or a
 $\bold{ruled\, manifold}$) are diffeomorphic. 
  These classification results on the symplectic structures on  $\mathbb{C}P^2$ or 
 $\bold{ruled\, manifolds}$  are due to Taubes and Lalonde-McDuff (\cite{T}, \cite{LMc}).

         In \cite{L} and in this paper, we assumed that  $\dim$ $H^2(M)<3$. Under this assumption, we defined $\bold{fixed\, point\, data}$.
   Notice that  the definition of $\bold{twist}$ does not depend on this restriction, it is defined when the fixed point set  
   consists only of surfaces.  
    As we mentioned in Remark~\ref{rem1'},  when there is an index 2 surface or the maximum is a surface in a certain ``type'', 
  the 1st Chern numbers of the normal bundles of these surfaces are uniquely determined by the fixed point data. The
    ``twist type'' controls the ``top'' and the ``bottom''.  
   If we remove the assumption $\dim$ $H^2(M)<3$, we have more fixed point set components. If the fixed point set consists only
  of surfaces, we still have our definition of  $\bold{twist}$. If the fixed point set also consists index 2 or index 4 isolated
  points,  then the reduced spaces change not only by diffeomorphism, but also by blow up or blow downs
    (\cite{GS}). The reduced spaces and their cohomologies are more complicated. The definition of 
 $\bold{twist}$ needs to be properly addressed.
  For the  definition of  $\bold{fixed\, point\, data}$, we
    need to specify more 1st Chern numbers of the normal bundles of the fixed point set components to determine 
    those  of the other fixed point set components. Alternatively, depending on the purposes, we may assume that all 
  these Chern numbers for all the fixed point
    set components  are compatibly given  and we include all of them  in the definition. For the construction of these more complicated 
  manifolds, we may not understand the symplectic structures on the reduced spaces well enough to do the local construction or to
 perform the global gluing. Certainly, in good cases, if our requrements are satisfied, we can still construct the manifolds with
  certain ``types'' other than those in  Theorem~\ref{thm1}.

         In summary, in order to construct a global manifold from a prescribed type, we need the ``data'' to be ``compatibly matched'', 
  and we need
   the reduced spaces to be ``good enough'' to allow us to perform the gluing. The idea itself can be applied to more general constructions of
   Hamiltonian $S^1$ manifolds, for instance, for non-semi-free Hamiltonian circle actions on other dimensional symplectic manifolds with
  other fixed point set components. Or more generally, depending on what is given,  the idea itself may be applied to the setting with no Hamiltonian circle actions.

    The main result of this paper is      
    
     \begin{theorem}\label{thm2}
   There exists a manifold for any of the prescribed fixed point data and twist type as in  Theorem~\ref{thm1}, except that for type (6a), the genuses
   $g$ and $g_1$ must satisfy  $1+g_1-2g>0$.
 % There exists a manifold in each of the six types in Theorem~\ref{thm1}.
 % In type (3) and type (6b), where $b_{min}$ can vary, there exists a manifold for each fixed $b_{min}$.
 % In type (6a), where  $b_{min}$, $g$ and $g_1$ can vary, there exists a manifold for each fixed $b_{min}$, $g$ and $g_1$ with $1+g_1-2g>0$.
    \end{theorem} 
  
   Recall that for type (6a), we computed the dual class $\eta$ of the index 2 surface in the reduced space at the critical level in \cite{L}: 
 $\eta=(1+g_1-2g)x+2y$
  for even $n$, and  $\eta=(2+g_1-2g)x+2y$ for odd $n$ (see Section 3 for the meaning of $x$ and $y$). 
   If  $1+g_1-2g\leq 0$, then these dual classes can not represent the connected symplectically embedded index 2 surface in the reduced spaces.
 The class $2y$ (for even $n$) or $x+2y$ (for odd $n$) represents non-connected symplectic submanifolds. The classes  $\eta=(1+g_1-2g)x+2y$ or  $\eta=(2+g_1-2g)x+2y$ with $1+g_1-2g<0$ can not represent the symplectically
 embedded image of the genus $g_1$ index 2 surface due to the following reason: the projection of the embedded surface (represented by  $\eta=(1+g_1-2g)x+2y$ or 
 $\eta=(2+g_1-2g)x+2y$)  to the ``base'' of the $S^2$ fibration (represented by
 $y$) is a smooth degree 2 map from a genus $g_1$ Riemann surface to a genus $g$ Riemann surface. By the analogue  Hurwitz formula for smooth
 maps between Riemann surfaces (see Theorem 1 in \cite{HY}), $g_1$ and $g$ must satisfy  $1+g_1-2g\geq 0$.  This is the reason we can not
  construct such manifolds.  This fact also helps us to strengthen Theorem~\ref{thm1}: 
  \footnote{The construction procedure gave me a good reason to study the symplectic submanifolds of the reduced spaces.
   This helped me to discover this piece of information.}
  \begin{corollary}
   For type (6a) in Theorem~\ref{thm1}, the genuses of the surfaces must satisfy  $1+g_1-2g>0$.
  \end{corollary}
 
      Theorem~\ref{thm1} and this corollary give us a complete list of the possible ``types''.  Theorem~\ref{thm2} shows that there exist 
 manifolds of all these ``types''.
   Together, they give us exactly the ``types'' of manifolds that can occur under the assumptions of  Theorem~\ref{thm1}.

      \subsubsection*{Acknowledgement}
    I would like to thank Denis Auroux for explaining to me the symplectic submanifolds of ruled surfaces. I thank him for explaining to me the
  Hurwitz formula for complex maps between Riemann surfaces. For type (6a), the idea of the proof  of the non-existence of the symplectic submanifolds
  in ruled surfaces for the cases $1+g_1-2g\leq0$ is due to him.    
    
    I would like to thank the referee for his invaluable comments and constructive suggestions.    
     
 \section{Construction of neighborhoods of critical and regular level sets of the moment map} 
    Construction of a neighborhood  of the minimum:

    If the minimum is isolated, by the equivariant Darboux-Weinstein theorem, a neighborhood  of the minimum is equivariantly symplectomorphic to $\mathbb{C}^3$ 
    with the standard symplectic form and $S^1$ action with weights $1, 1, 1$.

    If the minimum is a surface $\Sigma_g$, then a neighborhood of the minimum is a $\mathbb{C}^2$ bundle $E$ over $\Sigma_g$. It is
    easy to construct such a bundle with 1st Chern number $b_{min}$. (We can always construct 2 line bundles with Euler numbers $m$ and $n$ for any
    integer $m$ and $n$ and take the sum of the line bundles.). We can put a symplectic form on a neighborhood of the zero section of $E$ 
   (see \cite{GS}) and make $S^1$ act semi-freely by fixing
the base and rotating the fibers.
    More precisely, let $F$ be the unitary frame bundle of $E$ (assume there is a Hermitian metric on $E$), $T^*F$ be the cotangent bundle and $\omega_F$
    be the canonical symplectic form on $T^*F$. A connection on $F$ can give us the vertical cotangent bundle $V^*F$. 
   Let $i:V^*F\hookrightarrow T^*F$ 
    be the inclusion, and let $\omega_B$ be the symplectic form on the base $\Sigma_g$. Then $\omega_V=i^*\omega_F+\pi^*\omega_B$ is a 
    symplectic form on $V^*F$. The $U(2)$ action on
    $V^*F$ is Hamiltonian, so is the $U(2)$ action on $\mathbb{C}^2$. Take $V^*F^-\times\mathbb{C}^2$ with the product symplectic structure,
     where $V^*F^-$ is $V^*F$ with the form $-\omega_V$. There is a Hamiltonian $S^1$ action on this space: it acts trivially on the first factor and it acts
     on the second factor by rotation with weights $1, 1$. We can get a neighborhood of the zero section of $E$ with a symplectic form by reducing the $U(2)$ action at $0$. The semi-free 
     Hamiltonian $S^1$ action remains on $E$. By the equivariant Darboux-Weinstein theorem, any $S^1$-invariant symplectic structure on the above
  neighborhood is 
     equivariantly isomorphic to
     the above one. \\

    Notice that the above construction works for other dimensions  and for non-semi-free $S^1$ actions.

   Hence, we have the following lemma for the construction of a neighborhood of the minimum.
      \begin{lemma}\label{lem2.1}
      Assume $F$ is an isolated point or is a closed surface, then we can construct uniquely up to isomorphism a(n) (open) Hamiltonian $S^1$ 
   symplectic manifold $M_{[min,\epsilon_0)}$ such that $F$ is the minimum of the moment map, and such that $F$  has any prescribed 1st 
   Chern number of its normal bundle in  $M_{[min,\epsilon_0)}$ when $F$ is a surface.   
     \end{lemma}

    To construct a neighborhood of a critical level which contains an index 2 point or an index 2 surface, we can use the result of 
  Guillemin and Sternberg (\cite{GS}):

    Assume we have $(M_0, \omega_0, P, X)$, where $M_0$ is a compact symplectic manifold with symplectic form $\omega_0$, $X$ is a symplectic submanifold of $M_0$, and $P$ is a principal
    $S^1$ bundle over $M_0$, $\pi: P\rightarrow M_0$. Given a small open interval $I=(-\epsilon,+\epsilon)$, we can construct a symplectic manifold $M_I$ with a semi-free
    Hamiltonian  $S^1$ action with proper moment map $\phi$ which maps $M_I$ onto $I$. The symplectic form on $M_I$ is determined
  by $\omega_0$, the normal bundle of $X$ in $M_0$ and the Euler class of the circle bundle $P$ over $M_0$. Moreover, $0$ is an index 2 critical level
 of  $\phi$. For $-\tau<0$, $\phi^{-1}(-\tau)$ is diffeomorphic as an $S^1$ space
    to $P$,  the reduced spaces of $M_I$ at $-\tau$ and $0$ are diffeomorphic to $M_0$, and the reduced space  of $M_I$ at $\tau>0$ is the blow up of $M_0$ along $X$. The reduced symplectic forms
    satisfy the Duistermaat-Heckman formula (\cite{DH}). Moreover, $M_I$ is unique up to $S^1$- equivariant symplectomorphisms 
 (Theorem 13.1 in \cite{GS}).

   We summarize this in 
  \begin{lemma}\label{lem2.2}
      Assume we have $(M_0, \omega_0, P, X)$, where $M_0$ is a compact symplectic manifold with symplectic form $\omega_0$, $X$ is a symplectic submanifold of $M_0$, and $P$ is a principal $S^1$ bundle over $M_0$, $\pi: P\rightarrow M_0$. Given a small open interval $I=(-\epsilon,+\epsilon)$,
   we can construct uniquely up to isomorphism a semi-free Hamiltonian  $S^1$  symplectic manifold $M_I$ with moment map $\phi$ (which
 maps  $M_I$ onto $I$)  such that $\phi^{-1}(0)$ is a critical
 level with an index 2 fixed point set which is diffeomorphic to $X$,  $\phi^{-1}(-\tau)$ is diffeomorphic as an $S^1$ space
    to $P$, and  the reduced spaces of $M_I$ at $-\tau$ and $0$ are diffeomorphic to $M_0$, the reduced space  of $M_I$ at $\tau>0$ is 
the blow up of $M_0$ along $X$.
  \end{lemma}
  
    We can similarly construct a neighborhood of a critical level set which contains an index 4 fixed point
set component by looking at it upside down.\\

    The construction of a neighborhood  of the maximum is similar to the construction of a neighborhood  of the minimum.\\
    
      The above construction of neighborhoods of critical levels can only be done for small intervals of 
 critical values of $\phi$.   When we construct our manifolds, we may need a ``longer''
 neighborhood of a regular level. For this purpose, we construct a ``regular neighborhood''.

    If $I$ is a connected open interval consisting of regular values of $\phi$, then $\phi^{-1}(I)$ is $S^1$-equivariantly diffeomorphic to
    $\phi^{-1}(a)\times I$ where $a\in I$. The $S^1$-invariant symplectic form on $\phi^{-1}(I)$ is determined by the  symplectic 
    forms on the reduced spaces and the Euler class of the principal circle bundle $\phi^{-1}(a)$ over $M_a$. 
    \begin{lemma}(\cite{Mc1})    For an interval $I$ consisting of regular values of $\phi$, if $\omega$ is the symplectic form on $\phi^{-1}(I)$, and $\tau_s$ is the reduced symplectic form at
    $s\in I$, then
 $\omega$ is determined by the forms $\tau_s$ up to an $S^1$ equivariant diffeomorphism which preserves the level sets of 
    $\phi$.
    \end{lemma}
    
    See \cite{Mc1} for details.

     Assume that these neighborhoods we constructed above are
    ``big enough" so that we can shrink them a little such that they are manifolds with boundaries, the boundaries being regular level sets 
    of $\phi$.
    We will call these neighborhoods with boundaries of critical levels and regular levels ``local pieces''.
\section{Problems for local and global constructions} 
   
  We have seen in Section 2 what we need to prescribe to do the local construction. Let us review in this section what information the ``types'' in
   Theorem~\ref{thm1} can give us and what we need to take care of when we do the local and global constructions.

  Assume $a$ is a regular value of the moment map $\phi$, then $S^1$ acts freely on $\phi^{-1}(a)$, so we have a principal circle bundle

\begin{equation}\label{eq1}    
\begin{array}{ccl}
        S^1 & \hookrightarrow & \phi^{-1}(a)=P_a    \\                           
            &  & \downarrow    \\          
                           & &   \phi^{-1}(a)/S^1=M_a, \\
 \end{array}                          
  \end{equation}                             
where $M_a$ is the reduced space at $a$.
  
  Assume the minimum (or the maximum) of the moment map is a surface $\Sigma_g$. Let $a$ be a regular value right above the minimal value (or
  right below the maximal value) such that there are no other critical sets in the neighborhood except the minimum (or the maximum).
  Then $M_a$ is the total space of a sphere bundle over $\Sigma_g$. Let $S^2\times\Sigma_g$ be the trivial bundle, and let $E_{\Sigma_g}$ be
  the non-trivial bundle. Let $e(P_a)$ be the Euler class of the circle bundle $P_a$ over $M_a$.

     For $S^2\times\Sigma_g$, let $x, y\in H^2(S^2\times\Sigma_g)$ be a basis which consists of the dual classes of the fiber $S^2$ 
     and the base $\Sigma_g$ respectively. Then $\int_{base}x=\int_{S^2\times\Sigma_g}xy=1, \quad  \int_{fiber}y=\int_{S^2\times\Sigma_g}yx=1, \quad 
      \int_{fiber}x=\int_{S^2\times\Sigma_g}x^2=0, \hbox{and}  \int_{base}y=\int_{S^2\times\Sigma_g}y^2=0$.

       Similarly, for the non-trivial bundle $E_{\Sigma_g}$, let $x, y\in H^2(E_{\Sigma_g})$ be a basis which consists of the dual classes of the fiber and the 
      section $\Sigma_-$ which has self-intersection $-1$ respectively. Then $\int_{\Sigma_-}x=\int_{E_{\Sigma_g}}xy=1, \quad  \int_{fiber}y=\int_{E_{\Sigma_g}}xy=1, \quad 
       \int_{fiber}x=\int_{E_{\Sigma_g}}x^2=0,  \hbox{and}  \int_{\Sigma_-}y=\int_{E_{\Sigma_g}}y^2=-1$.

     We have the following lemmas
    \begin{lemma} \label{lem1} (\cite{L})
       If the minimum is a surface $\Sigma_g$ with 1st Chern number $b_{min}$ of its normal bundle, then $M_a$ (a above the minimum) is diffeomorphic 
       to $S^2\times\Sigma_g$ if and only if $b_{min}=2k$ is even, and it is diffeomorphic to $E_{\Sigma_g}$ if and only if $b_{min}=2k+1$ is odd.
       In either case, $e(P_a)=kx-y$.
       \end{lemma}
 
  \begin{lemma}\label{lem2} (\cite{L})
     If the maximum is a surface $\Sigma_g$ with 1st Chern number $b_{max}$ of its normal bundle, then $M_a$ (a below the maximum) is diffeomorphic to 
     $S^2\times\Sigma_g$ if and only if $b_{max}=2k'$ is even, and it is diffeomorphic to $E_{\Sigma_g}$ if and only if 
     $b_{max}=2k'+1$ is odd. In either case, $e(P_a)=-k'x+y$.\\ 
     \end{lemma}
 Let $c$ be a critical value of the moment map $\phi$, and  $S\subset \phi^{-1}(c)$ be an index 2 fixed surface.  
 Assume for small $\epsilon$, $S$ is the only critical set in $\phi^{-1}([c-\epsilon, c+\epsilon])$. 
 Let $P_{c-\epsilon}$ be the circle bundle $\phi^{-1}(c-\epsilon)$ over $M_{c-\epsilon}$, $e(P_{c-\epsilon})$ be its Euler class.
 $P_{c+\epsilon}$ and $e(P_{c+\epsilon})$ have similar meanings. By \cite{GS} or \cite{Mc1}, 
 $M_{c+\epsilon}$ is diffeomorphic to $M_{c-\epsilon}$ (in dimension 6), let's call them $M_{red}$ (reduced space).
  Moreover, there is a symplectic embedding
$i: S\rightarrow M_{red}$ such that  if $i(S)=Z$, and $\eta$ is the dual class of $Z$ in 
             $M_{red}$, then
 \begin{equation} \label{eq2}           
             e(P_{c+\epsilon})=e(P_{c-\epsilon})+\eta.                  
 \end{equation}
 
    By \cite{GS}, when the moment map crosses an index 2 point (index 4 point), the reduced space changes by a blow up at a point (blow down of an exceptional
    divisor). Now let $M_-$ be the reduced space before surgery (blow up or blow down), $M_+$ be the reduced space after surgery, and let
    $e(P_-)$ be the Euler class of the principal circle bundle over $M_-$, $e(P_+)$ be the Euler class of the principal circle bundle over $M_+$. 
    Let $\eta_{ex,div}$ be the dual class of the exceptional sphere in $M_+$ after blow up or in $M_-$ before blow down. Let's always call the blow
     down map by $\beta$. By \cite{GS}, we have the following:
     
      Corresponding to blow up,
    \begin{equation} \label{eq3}
       e(P_+)=\beta^*e(P_-)+\eta_{ex,div}.
    \end{equation}

      Corresponding to blow down,
    \begin{equation} \label{eq4}
       \beta^*e(P_+)=e(P_-)+\eta_{ex,div}.
    \end{equation}

For each type of manifolds in  Theorem~\ref{thm1}, we can use Lemma~\ref{lem1}, Lemma~\ref{lem2} and equations (\ref{eq2}), (\ref{eq3}), (\ref{eq4}) to compute  
 the Euler classes of the principal $S^1$ bundles (\ref{eq1})
 for regular levels in terms of the fixed point data.  By Duistermaat-Heckman formula (\cite{DH}), we can compute the cohomology classes of the reduced symplectic forms on the reduced spaces.\\

       Recall that for the types of manifolds in Theorem~\ref{thm1}, the reduced spaces
are all smooth manifolds (see Section 10 in \cite{GS}). Moreover, from above we see that (except at the minimal value and at the maximal value)
   they are diffeomorphic to either
$\mathbb{C}P^2$ or the total space of an $S^2$ bundle over a Riemann surface (assume for manifold of type (5) the moment map crosses the index
4 point before it crosses the index 2 point). Now, let us list the problems we should take care of for the construction:\\
 $\mathbf{1.}$ The  cohomology classes of the reduced forms given by Duistermaat-Heckman formula are in terms of the fixed
point data, the size of the minimum (maximum) if the minimum (maximum) is a 
surface, and the distances between the moment map images of the fixed point set components. The question is which cohomology class
 in the reduced space contains a symplectic form.\\
$\mathbf{2.}$  If $F$ is an index 2 fixed surface, then $F$ is symplectically embedded in
the reduced space at the critical level (and nearby regular levels). As we have seen (see \cite{L} for details), in each case, 
the dual cohomology class of $F$ in the reduced space can be computed in terms of the fixed point data. The question is  which dual class 
can be represented by a symplectic submanifold.\\
$\mathbf{3.}$ When we glue the ``local pieces'' together, we glue
them along a regular level. Since this regular level is a principal $S^1$
bundle over the reduced space, we can glue the two reduced spaces on the 
local pieces together (the fibers will be identified correspondingly).
We only know the cohomology classes of the reduced symplectic form are
the same on these two reduced spaces, the symplectic forms may be chosen differently.
The question is if there is a diffeomorphism between the two reduced spaces
such that one form is pulled to another by this diffeomorphism. 

    Section 4 and Section 5 will give answers to the above questions.
   \section{Symplectic $\mathbb{C}P^2$, symplectic ruled manifolds, and
 their symplectic submanifolds}
   \begin{lemma}\label{lem9.3.1}
   Let $f: (M, \omega_1)\rightarrow (M, \omega_2)$ be a diffeomorphism
  such that $f^*\omega_2=\omega_1$. If $A\subset (M, \omega_1)$ is a symplectic submanifold,
  then $f(A)\subset (M, \omega_2)$ is a symplectic submanifold.
 \end{lemma}
  \begin{proof}
  $f_*(TA)=T(f(A))$ is an isomorphism. Take $X, Y\in TA$,
 then $\omega_2(f_*X, f_*Y)=f^*\omega_2(X, Y)=\omega_1(X, Y)$.
 \end{proof} 

  $\mathbb{C}P^2$ with its standard complex structure carries  K$\ddot{a}$hler forms. Each complex submanifold of   $\mathbb{C}P^2$ is symplectic with
 respect to these  K$\ddot{a}$hler forms.
 
   For the symplectic structures on  $\mathbb{C}P^2$, we have
    \begin{theorem}\label{T}
    (\cite{T}) If $(M,\omega)$ is a symplectic manifold and $M$ is smoothly diffeomorphic to $\mathbb{C}P^2$, then $\omega$ is diffeomorphic to a standard
    K$\Ddot{a}$hler form. 
    \end{theorem}

 \begin{lemma}\label{lem9.3.2}
         
   Let $u\in H^2(\mathbb{C}P^2)$ be a generator.
Then a cohomology class $a=lu$ contains a symplectic form of  $\mathbb{C}P^2$ 
    if and only if $l\neq 0$.
 \end{lemma}
  \begin{proof}
  Any non-zero multiple of the standard symplectic form on  $\mathbb{C}P^2$ is symplectic.
 \end{proof}
 \begin{lemma}\label{lem9.3.3}
   For $(\mathbb{C}P^2,\omega)$, let $u\in H^2(\mathbb{C}P^2, \mathbb{Z})$ be a generator.  Then any class $\eta=au$ with $a$ being a positive integer can be represented by a connected
       symplectic submanifold.
 \end{lemma}
 \begin{proof}
 By Theorem~\ref{T}, $\omega$ is diffeomorphic to a stardard K$\Ddot{a}$hler form. Since  $\eta=au$ with $a$ being a positive integer represents a connected
complex curve for the standard complex structure,  it represents a connected
symplectic submanifold for the standard K$\ddot{a}$hler  form.  The lemma
follows from Lemma~\ref{lem9.3.1}. 
\end{proof}

 \begin{definition}\label{def9.3.1}
       An $S^2$ fibration over a Riemann surface $\pi: E\rightarrow \Sigma$ is called a $\bold{ruled\, manifold}$.
 \end{definition}
   Ruled manifolds are symplectic manifolds. One can choose the fibration such that the symplectic forms are nondegenerate on the fibers.    
 \begin{definition}\label{def9.3.2}
       A symplectic form on a ruled manifold is $\bold{compatible}$ with the ruling $\pi: E\rightarrow \Sigma$ if it is nondegenerate on the fibers. 
  \end{definition}
  The following fact is well known
   \begin{lemma}\label{lem9.3.4}
    Let $E$ be a ruled manifold, and  $F$ be the fiber of $E$.  Any cohomology class $a\in H^2(E)$
    which satisfies the following conditions carries a symplectic form 
    which is nondegenerate on the fibers.\\
    (a) $a(F)$ and $a^2(E)$ are positive, and\\
    (b) $a^2(E)>(a(F))^2$ if the bundle is non-trivial.
  \end{lemma}
  
  With this symplectic form and a compatible complex structure, 
every connected complex submanifold is a connected symplectic submanifold.

   By D. McDuff, we have 
    \begin{lemma}\label{lem9.3.4'}
  (\cite{Mc2})  
   Let $E$ be a ruled manifold, and  $F$ be the fiber of $E$. Then the cohomology class $a\in H^2(E)$ of any
   symplectic form on $E$ which is compatible with the ruling 
    satisfies the  conditions\\
    (a) $a(F)$ and $a^2(E)$ are positive, and\\
    (b) $a^2(E)>(a(F))^2$ if the bundle is non-trivial.
  \end{lemma}

  For the symplectic forms in the same cohomology class on a ruled manifold,
 we have  

  \begin{theorem}\label{LMc}
     (\cite{LMc}) Let $\omega_0, \omega_1$ be two cohomologous symplectic forms on the ruled 4-manifold $\pi: E\rightarrow \Sigma$, then there is a diffeomorphism
    $\Phi$ of $E$ such that $\Phi^*(\omega_0)=\omega_1$. Moreover, if we assume that the forms $\omega_0, \omega_1$ are both compatible with $\pi$, 
    then they are isotopic.
    \end{theorem}

   We have the following lemmas about the connected symplectic submanifolds of ruled surfaces.
     \begin{lemma}\label{lem9.3.5}
     For $(\Sigma_g\times S^2,\omega)$, where $\omega$ is compatible with the ruling, 
  let $x, y$ be the dual class of the fiber $S^2$ and the dual class of the base $\Sigma_g$ respectively. 
     Then $x$, $y$,
         and any class $ax+by$ with $a, b$ being positive integers can be represented by a connected symplectic submanifold.

            If $\int{[\omega]\cdot y}>n \int{[\omega]\cdot x}$ ($n=1,2,...$),
          then $y-nx$  can be represented by a connected symplectic submanifold.
    \end{lemma}
    
    \begin{lemma}\label{lem9.3.6} 
      For the non-trivial $S^2$ bundle $(E_{\Sigma_g},\omega)$ over $\Sigma_g$, where $\omega$ is compatible with the ruling,
 let $x, y$ be the dual class of the fiber $S^2$ and the dual class of
        the section with self-intersection $-1$ respectively. Then $x$, $y$, $x+y$, and any class $ax+by$ with integers $a>1, b>0$ can be 
        represented by a connected
        symplectic submanifold.

        If $\int{[\omega]\cdot y}>n \int{[\omega]\cdot x}$ ($n=1,2,...$), then $y-nx$ 
          can be represented by a connected symplectic submanifold.     
     \end{lemma}
      
     Proof of Lemma~\ref{lem9.3.5}:
     \begin{proof}
     The classes  $x$ and $y$ clearly represent symplectic submanifolds.  For the classes $ax+by$ with $a, b$ being positive integers,
    we can construct their symplectic representatives. Take the class $p_i\times S^2$, where $i=1, 2, ..., a$, and take the
    class $\Sigma_g\times q_j$, where $j=1, 2, ..., b$. 
    These are symplectic submanifolds and they intersect
    positively transversely at $ab$ points. Then one can smooth out the intersection points to obtain a smooth submanifold which
    is symplectic and represented by the class $ax+by$.
      
         When the symplectic area of the base is $n$ times larger than that of the fiber, the class $y-nx$ represents the 
     section at $\infty$ of the projectivised bundle $P(L\oplus\mathbb{C})$, where $L$ is a complex line bundle of degree $2n$
     over $\Sigma_g$. So it represents a symplectic submanifold.   
     \end{proof}
  
    Alternatively, one may choose compatible complex structures on ruled manifolds and argue that some classes represent 
   complex submanifolds.   
    The proof  of Lemma~\ref{lem9.3.6} is similar.

   \section{Gluing neighborhoods together}
   
   In Section 2, we constructed ``local pieces". 
  We will glue these ``local pieces" together along their boundaries.
   Let $L_1$ and $L_2$ be two ``adjacent'' ``local pieces'' such that $\phi|L_1\leq\phi|L_2$.  Let $a$ be the regular
  value such that  $\phi^{-1}(a)$ is the level along which we glue $L_1$ and $L_2$. 
  Suppose $\phi^{-1}_-(a)$ is the boundary of $L_1$, 
  $\phi^{-1}_+(a)$ is the boundary of $L_2$, and we will identify them
  to get  $\phi^{-1}(a)$. Now let $M_a^-=\phi^{-1}_-(a)/S^1$, and 
  $M_a^+=\phi^{-1}_+(a)/S^1$. Then  $\phi^{-1}_-(a)$ is a principal circle bundle over $M_a^-$, and $\phi^{-1}_+(a)$ is a principal circle
  bundle over  $M_a^+$. Let  $\omega_-$  and $\omega_+$ be the reduced symplectic forms on  $M_a^-$ and on $M_a^+$ respectively. 
  We state the following lemma for convenience, it is obviously valid. 
 \begin{lemma}\label{lem5.1}
  Assume that the two principal circle bundles  $\phi^{-1}_-(a)\rightarrow M_a^-$ and $\phi^{-1}_+(a)\rightarrow M_a^+$ have the same
  Euler class. If there is a symplectomorphism $g$ between $(M_a^-, \omega_-)$  and   $(M_a^+, \omega_+)$, then $g$ induces an 
 isomorphism between  $\phi^{-1}_-(a)$ and  $\phi^{-1}_+(a)$.
 \end{lemma}
   
   In our case, of course the two circle bundles have the same Euler class.
   As we have seen in Section 3, our reduced spaces are diffeomorphic to $\mathbb{C}P^2$ or diffeomorphic
  to a ruled manifold.   $M_a^-$ and $M_a^+$ are diffeomorphic as
  smooth manifolds.   The cohomology classes of the reduced
  symplectic forms $\omega_-$ on $M_a^-$ and $\omega_+$ on $M_a^+$ agree. 
  By Theorem~\ref{T} or by Theorem~\ref{LMc}, there is a 
  diffeomorphism $g: (M_a^-, \omega_-)\rightarrow (M_a^+, \omega_+)$ such
  that $g^*\omega_+=\omega_-$.
   By Lemma~\ref{lem5.1},  this symplectomorphism $g$ between the reduced
  spaces induces an isomorphism between  $\phi^{-1}_-(a)$ and  $\phi^{-1}_+(a)$. We call the identified level set  $\phi^{-1}(a)$.

 \section{Proof of Theorem~\ref{thm2}}

  We are now ready to construct our manifolds.  It will boil down to the question of choosing some parameters suitably.

   \begin{proof}
     In the following, we always use $\alpha_0$ to denote the size of the minimum if the minimum is a surface, and use $t_0, t_1, t_2$ to denote the distances
between the moment map images of the fixed point set components as the moment
map goes from the minimum to the maximum. For convenience, we will use the following pictures. On the
    left, between the critical sets, we indicate the Euler class of the principal $S^1$ bundle over the reduced space at that regular interval, and on the
    right of the index 2 surface, we indicate the dual class of the index 2 surface in the reduced space.\\

    By Section 5, the gluing of the ``local pieces'' can always be done, therefore
 we won't ``glue explicitly'' in the following proof, we will check requirements $\mathbf{1.}$ and $\mathbf{2.}$ listed in Section 3. \\

   For manifold of type (1), we have seen an example in \cite{L}. Let us use the construction method to prove its existence. 
  We have the following picture ($u\in H^2(\mathbb{C}P^2)$ is a generator):

  \begin{figure}[h!]
    \scalebox{.60}{\includegraphics{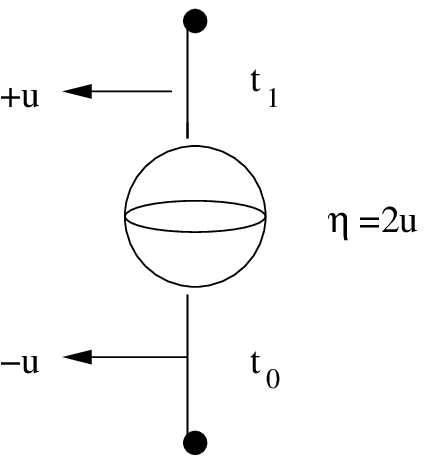}}
    \end{figure}
The reduced space at any non-extremal level is $\mathbb{C}P^2$.

 Since the minimum and the maximum are unique,  $t_0>0$, and $t_1>0$.

  Let us  use Duistermaat-Heckman formula (\cite{DH}) to compute the cohomology classes of
the reduced symplectic forms. If $0<t\leq t_0$, then $[\omega_t]=ut$, by Lemma~\ref{lem9.3.2}, for each t
in this interval, $[\omega_t]$ contains a symplectic form.
 
   By Lemma~\ref{lem9.3.3}, the index 2 sphere is a symplectic submanifold of
  $\mathbb{C}P^2$.

   Starting new at $t_0$, by  Duistermaat-Heckman formula, for $0<t\leq t_1$, $[\omega_t]=t_0u-ut$. 
 Since $[\omega_{t_1}]=0$,  $t_1=t_0$. By lemma~\ref{lem9.3.2} again, each $[\omega_t]$ in this interval contains a symplectic form.

   This proves the existence of  manifold of type (1) in Theorem~\ref{thm1}.  
    
   We can similarly prove the existence of  manifold of type (2) in  Theorem~\ref{thm1}.\\

 For manifold(s) of type (3),  we saw
in \cite{L} that only when $b_{min}=1$, $\phi$ can cross the index 4 point before it crosses the index 2 sphere. 
  In this case, we constructed a toric variety. So let us assume the moment map crosses the critical sets in the
  following order.
  
      \begin{figure}[h!]
    \scalebox{.60}{\includegraphics{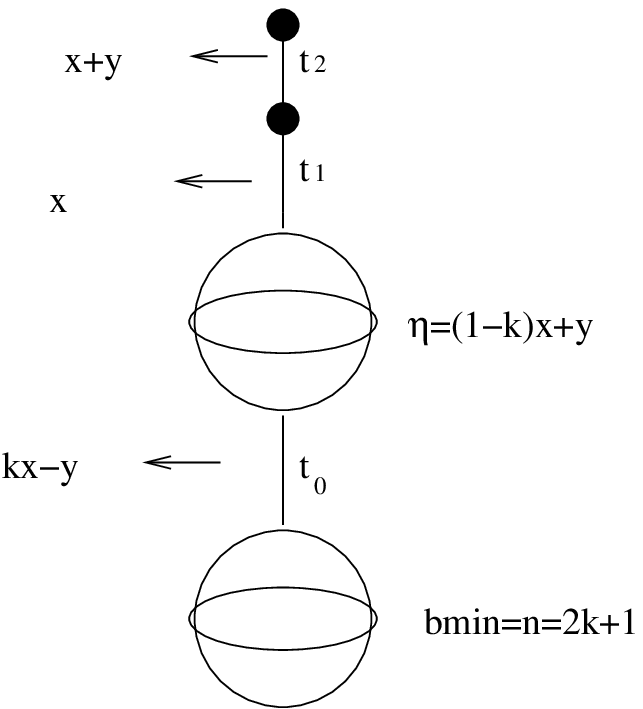}}
    \end{figure}

  Since the minimum and the maximum are unique,  $t_0>0, t_2>0$.

    Note that up to the index 4 critical level (excluding this critical level and the minimal level), the reduced space is $\tilde{\mathbb{C}P^2}$ ($\mathbb{C}P^2$
 blown up at a point). 
  By Duistermaat-Heckman formula, 
   if $0<t\leq t_0$, then $[\omega_t]=\alpha_0x-(kx-y)t=(\alpha_0-kt)x+ty$. To make $[\omega_t]$ contain symplectic forms for each $t$, by Lemma~\ref{lem9.3.4}, we can
     choose $\alpha_0>(1+k)t$ for each $t$ in this interval. Since $\alpha_0>0$, we only need $\alpha_0>(1+k)t_0$.

      By Lemma~\ref{lem9.3.6}, if $k\leq 1$, then $\eta$ represents a connected symplectic submanifold (this imposes no condition);
      if $k>1$, then if $\alpha_0-2kt_0>0$,  $\eta$  represents a connected symplectic submanifold.

     For the next interval, let's start with $t=0$. So when $0\leq t<t_1$ (the
reduced space $M_{red}$ at the index 4 critical level is $\mathbb{C}P^2$.),
      $$ [\omega_t]=[\omega_{t_0}]-xt=(\alpha_0-kt_0-t)x+t_0y.$$
     To make $[\omega_t]$ contain a symplectic form for each $t$,  we  choose $\alpha_0-(k+1)t_0>t$ for $0\leq t<t_1$.

     After $\phi$ crosses the index 4 point, $M_{red}$ is $\mathbb{C}P^2$ and $x+y$ is a generator of $H^2(\mathbb{C}P^2)$.
      If we do the computation starting from the top, we see that each
     cohomology class in this interval contains a symplectic form. To make the data compatible, we need
     $$(\alpha_0-kt_0-t_1)x+t_0y=t_2(x+y),$$
     i.e.,  $ t_2=t_0$,  and $\alpha_0=t_1+(k+1)t_0$.

     In summary, we choose the parameters such that 
     $$\alpha_0>2kt_0, \quad \alpha_0=t_1+(k+1)t_0, \quad t_1>0  \quad (t_0=t_2).$$
     For each fixed $k$, we can adjust $\alpha_0, t_0, t_1$ such that the above condition is satisfied
     This proves the existence of manifold(s) of type (3).\\

     We have shown in \cite{L} that manifold of type (4) is diffeomorphic to $\mathbb{C}P^3$.\\
     
     For manifold of type (5), let us prove the existence for the case when $\phi$ crosses the index 4 point before it crosses the index 2 point.
 
           \begin{figure}[h!]
    \scalebox{.60}{\includegraphics{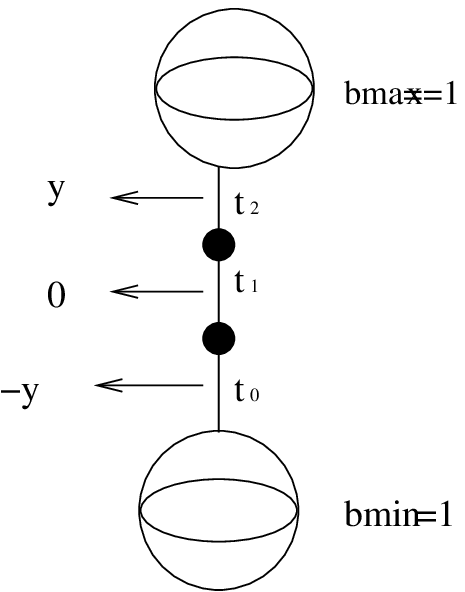}}
    \end{figure}

     If $0<t<t_0$, $M_{red}=\tilde{\mathbb{C}P^2}$, then $[\omega_t]=\alpha_0x+ty$, by Lemma~\ref{lem9.3.4}, we may choose  
     $\alpha_0>t$ for each $t$ in this interval to make $[\omega_t]$ contain
symplectic forms. Since $M_{red}$ at the index 4 point critical level is $\mathbb{C}P^2$ ($x+y\in H^2(\mathbb{C}P^2)$ is a generator), to
     make the cohomology class contain a symplectic form for $t=t_0$,
by Lemma~\ref{lem9.3.2},  we choose $\alpha_0=t_0$.

      Notice that the Euler class of the principal $S^1$ bundle over $M_{red}$ for the $t_1$ interval is $0$. We can do
       a similar computation starting from the top and choose $\tilde\alpha=t_2$, where $\tilde\alpha$ is the size of the
       maximum. To make the data compatible, (the cohomology class of the reduced symplectic form  at the index 4 critical level and at the index 2 critical level
       should be the same) we should choose $t_0=t_2$.

       In summary, we can choose parameters such that 
       $$\alpha_0=\tilde\alpha=t_0=t_2.$$
       This proves the existence of manifold of type (5).\\

       For manifolds of type (6a), we will prove the existence for the case when $b_{min}=n=2k$ is even
       (the reduced space is the trivial $S^2$ bundle over $\Sigma_g$ at any non-extrem level.). The proof for the case when $b_{min}$ is odd is similar (the
       reduced space is the non-trivial $S^2$ bundle over $\Sigma_g$ at any non-extrem level).

                  \begin{figure}[h!]
    \scalebox{.60}{\includegraphics{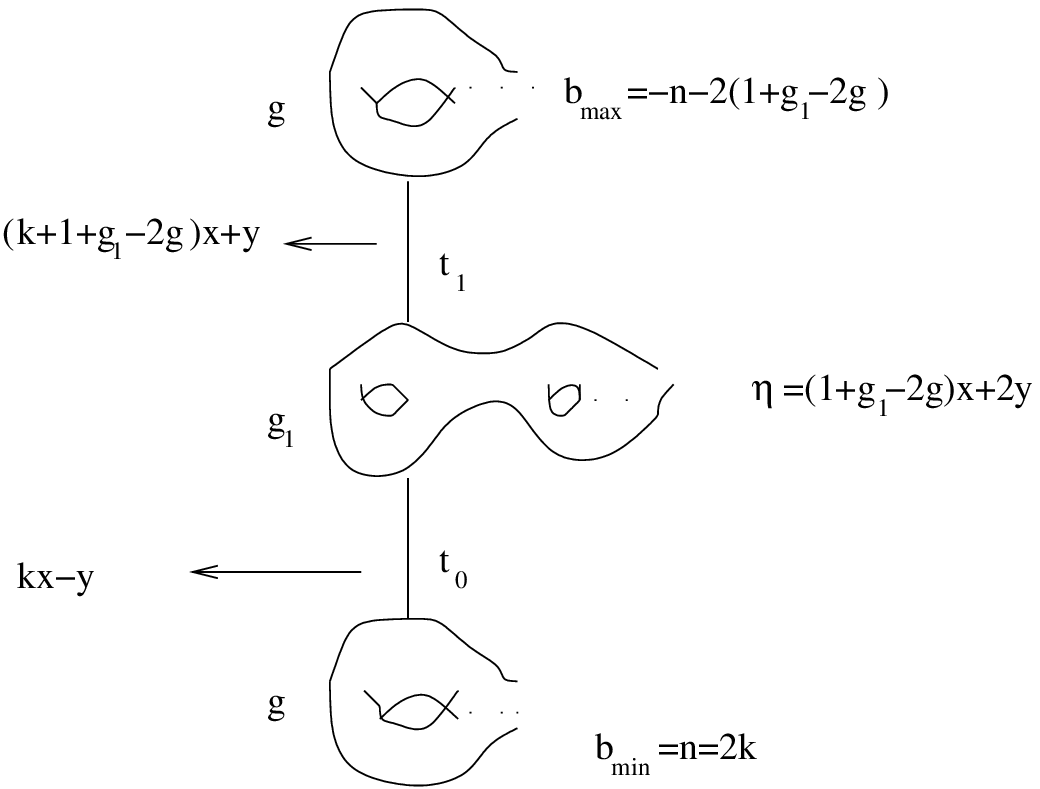}}
    \end{figure}

       For $0<t\leq t_0$, $[\omega_t]=\alpha_0x-(kx-y)t=(\alpha_0-kt)x+ty$.
By Lemma~\ref{lem9.3.4}, we need $\alpha_0-kt>0$ for each $t$ to make each cohomology
       class contain a symplectic form. Since $\alpha_0>0$, we only need $\alpha_0>kt_0$.

       By Lemma~\ref{lem9.3.5}, when $1+g_1-2g>0$, the index 2 surface is a symplectic submanifold of the reduced space.

       For $0<t\leq t_1$, $[\omega_t]=(\alpha_0-kt_0)x+t_0y-[(k+1+g_1-2g)x+y]t=[\alpha_0-kt_0-(k+1+g_1-2g)t]x+(t_0-t)y$. We need
       $t_0=t_1$, and $\alpha_0-kt_0-(k+1+g_1-2g)t_1>0$, i.e., $\alpha_0-nt_0>(1+g_1-2g)t_0$.
        If $k\geq 0$, it is easy to see this choice makes
       each cohomology class in this interval contain a symplectic form. If $k<0$, then we need $\alpha_0-kt_0-kt>(1+g_1-2g)t$
       for each $0<t<t_0$ to make the cohomology classes contain symplectic forms.

       In summary, we make the following choices of the parameters:
       $$\alpha_0-kt_0-kt>(1+g_1-2g)t>0,  t\in [0, t_0].$$
       We can adjust $\alpha_0$ and $t_0$ such that the above condition is satisfied for any fixed $k$ and $g, g_1$ such that
       $1+g_1-2g>0$. Therefore  manifolds of type (6a) exist when  $1+g_1-2g>0$.\\

       For manifolds (6b) (twisted case), $M_{red}$ is $S^2\times S^2$.

      When $b_{min}=n=2k$ and $b_{max}=0$, we have the following  picture:

             \begin{figure}[h!]
    \scalebox{.60}{\includegraphics{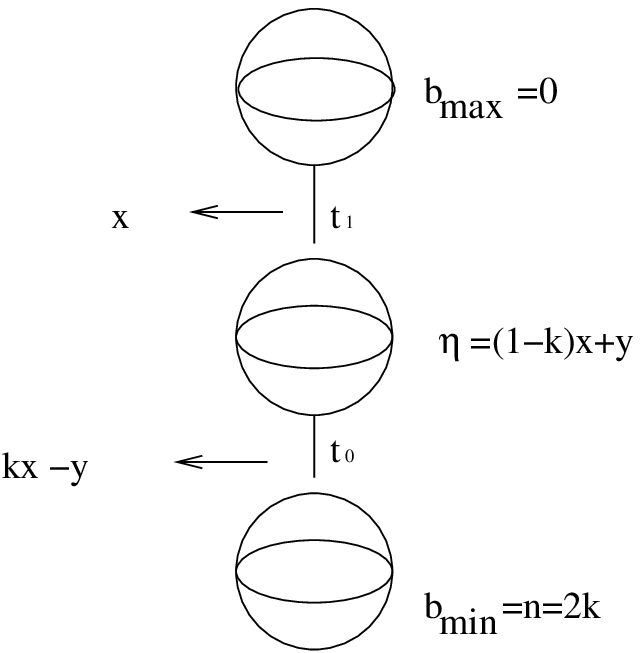}}
    \end{figure}

       Notice that $t_0>0$, and  $t_1>0$.

      When $0<t\leq t_0$, $[\omega_t]=\alpha_0x-(kx-y)t=(\alpha_0-kt)x+ty$, it is easy to see that when $\alpha_0>kt_0$, each $[\omega_t]$ contains
       a symplectic form.

       By Lemma~\ref{lem9.3.5}, if $k\leq 1$, then $\eta$ can be represented by a symplectic submanifold;
       if $k>1$, then when $\alpha_0-(n-1)t_0>0$, $\eta$ can be represented by a symplectic submanifold.

      When $0\leq t\leq t_1$ (new start at $t=0$), $[\omega_t]=(\alpha_0-kt_0)x+t_0y-xt=(\alpha_0-kt_0-t)x+t_0y$.
      Since there is a twist, when $t=t_1$, we should have  
      $$\alpha_0-kt_0-t_1=0.$$
      Notice that this kind of choice makes $[\omega_t]$ contain a symplectic form for each  $0<t<t_1$.

      In summary, we can choose 
      $$\alpha_0>(n-1)t_0, \quad t_1=\alpha_0-kt_0>0.$$
      For each fixed $k$, the choice is achievable. This proves the existence of manifold of type (6b).\\

      If this manifold exists, then its upside down also exists.

  \end{proof}

 \end{document}